\documentclass[10pt]{article}
\usepackage{graphicx}
\usepackage{amssymb}
\usepackage{epstopdf}
\usepackage{amsmath}
\usepackage{amsfonts}

\def\phi{{\varphi}}

\DeclareSymbolFont{AMSb}{U}{msb}{m}{n}
\DeclareMathSymbol{\N}{\mathbin}{AMSb}{"4E}
\DeclareMathSymbol{\Z}{\mathbin}{AMSb}{"5A}
\DeclareMathSymbol{\R}{\mathbin}{AMSb}{"52}
\DeclareMathSymbol{\Q}{\mathbin}{AMSb}{"51}
\DeclareMathSymbol{\I}{\mathbin}{AMSb}{"49}
\DeclareMathSymbol{\C}{\mathbin}{AMSb}{"43}

\def\be{\begin{equation}}
\def\ee{\end{equation}}
\def\ber{\begin{eqnarray}}
\def\eer{\end{eqnarray}}

\def\beq{\begin{equation}}
\def\eeq{\end{equation}}

\begin{document}

\addtolength{\textheight}{0 cm} \addtolength{\hoffset}{0 cm}
\addtolength{\textwidth}{0 cm} \addtolength{\voffset}{0 cm}

\newenvironment{acknowledgement}{\noindent\textbf{Acknowledgement.}\em}{}

\setcounter{secnumdepth}{5}
 \newtheorem{proposition}{Proposition}[section]
\newtheorem{theorem}{Theorem}[section]
\newtheorem{lemma}[theorem]{Lemma}
\newtheorem{coro}[theorem]{Corollary}
\newtheorem{remark}[theorem]{Remark}
\newtheorem{claim}[theorem]{Claim}
\newtheorem{conj}[theorem]{Conjecture}
\newtheorem{definition}[theorem]{Definition}
\newtheorem{application}{Application}

\newtheorem{corollary}[theorem]{Corollary}

\title{ Invariance  properties of the  Monge-Kantorovich mass transport  problem}
\author{
Abbas Moameni
\thanks{Supported by a grant from the Natural Sciences and Engineering Research Council of Canada.}
\hspace{2mm}\\
{\it\small School of Mathematics and Statistics}\\
{\it\small Carleton University}\\
{\it\small Ottawa, ON, Canada K1S 5B6}\\
{\it\small  momeni@math.carleton.ca}\\
%
}
\maketitle

\begin{abstract}
We consider the multidimensional Monge-Kantorovich transport problem in an abstract setting.  Our main results state that if a  cost function and
marginal measures are invariant by  a family of  transformations, then a solution of the Kantorovich relaxation  problem  and  a solution of its dual can be chosen so that they are
 invariant
 under the same family of transformations.  This provides a new tool to  study and analyze  the support of optimal transport plans and consequently to scrutinize the Monge problem.
Birkhoff's Ergodic theorem  is an essential tool in our analysis.
\end{abstract}

\section{Introduction}
We consider the Monge-Kantorovich transport problem  for Borel probability measures $\mu_1, \mu_2,..., \mu_n$ on Polish spaces   $X_1, X_2,..., X_n.$
The cost function $c: X_1 \times X_2 \times... \times X_n  \to [0, \infty]$ is Borel measurable. $\Pi(\mu_1,...,\mu_n)$ is the set of Borel probability measures on
$X_1 \times X_2 \times... \times X_n$  which have $X_i$-marginal $\mu_i$ for each $i \in \{1,2,...,n\}.$
The transport cost associated to a transport plan
 $\pi \in \Pi(\mu_1,...,\mu_n)$ is given by
\[I_c(\pi)=\int_{X_1 \times X_2 \times... \times X_n} c(x_1,...,x_n) \, d \pi.\]
We consider  the Monge-Kantrovich transport problem,
\[  \inf \{ I_c(\pi); \pi \in \Pi(\mu_1,...,\mu_n)\}. \qquad \qquad (MK)        \]
 The dual formulation of $(MK)$ takes the following aspect \cite{B-L-S}.
\begin{theorem}
 Assume that $X_1, X_2,..., X_n$ are Polish spaces equipped with probability measures $\mu_1, \mu_2,..., \mu_n,$
that $c: X_1 \times X_2 \times... \times X_n  \to [0, \infty]$ is Borel measurable and $\otimes_{i=1}^n\mu_i$-a.e.
finite and that there exists a finite transport plan. Then there exist a Borel measurable dual maximizer $(\phi_1,\phi_2,...,\phi_n),$ i.e.
 functions $\phi_i: X_i \to [-\infty, \infty)$ satisfying $ c(x_1,...,x_n) \geq\sum_{i=1}^n\phi_i(x_i)$ for all $(x_1,..., x_n) \in X_1 \times X_2 \times... \times X_n $
such that
\[\inf_{\gamma \in \Pi(\mu_1,..., \mu_n)} \int_{X_1 \times X_2 \times... \times X_n} c(x_1,...,x_n) \, d \gamma=\sum_{i=1}^n\int_{X_i}\phi_i(x_i) \, d \mu_i \]
\end{theorem}
Indeed, the authors in \cite{B-L-S} proved the latter theorem  for the two marginal case. The multi marginal case follows by the same argument 
and since it does not require  new ideas we do not elaborate. 
By virtue of the above theorem,  $(MK)$ is dual to the following problem.

\[ \sup\big \{ \sum_{i=1}^n\int_{X_i} \phi_i(x_i) \, d \mu_i ; (\phi_1,...,\phi_n) \in K_c \big \},    \qquad \qquad   (DK)\]
where $K_c$ is the set of  $(\phi_1,\phi_2,...,\phi_n),$ such that
 functions $\phi_i: X_i \to [-\infty, \infty)$ are Borel measurable and  $ c(x_1,...,x_n) \geq\sum_{i=1}^n\phi_i(x_i)$ for
all $(x_1,..., x_n) \in X_1 \times X_2 \times... \times X_n.$\\

Let $(X, \Sigma,\mu) $  be  a measure space.   A map   $T: X \to X$ is said to be a  $\mu$-measure  preserving  transformation if
 $T\# \mu=\mu,$ i.e.,
\[\forall f \in {\mathcal L}^1(\mu), \quad \, f \circ T \in {\mathcal L}^1(\mu) \, \& \, \int_{X}f(x) \, d\mu=\int_{X}f(Tx) \, d\mu,\]
where  ${\mathcal L}^1(\mu)$ is the set of integrable functions on $X.$ We also have the following definition for periodic maps.
\begin{definition}
  A map $T: X \to X$ is  called periodic of order $m$ if  $T^m$,  the $m$-th iterate of $T$, is the identity map  i.e.  $T^m(x)=x$  for every   $x \in X.$  We also say that
a map $T: X \to X$ is periodic if there  exists a positive integer $m$ such that $T$ is periodic of order $m$. The smallest such value of $m$  is called the period of $T.$
\end{definition}

If the    cost function  and  the
marginal measures  are invariant by  a family of measure preserving   transformations, then one expect the optimal mass
transport problem $(MK)$ and its dual $(DK)$ to possess  solutions  that
 are  invariant under the same transformations. Let us first state  our main results  for the invariance properties of $(DK).$
\begin{theorem}\label{main}  Let  $X_1, X_2,..., X_n$ be  Polish spaces equipped with probability measures $\mu_1, \mu_2,..., \mu_n,$
and  $c: X_1 \times X_2 \times... \times X_n  \to [0, \infty]$ be  Borel measurable and $\otimes_{i=1}^n\mu_i$-a.e.
finite. We also assume  that there exists a finite transport plan,  and  that  the dual problem $(DK)$ has  a solution $(\phi_1,...,\phi_n)$
such that $\phi_j \in {\mathcal L}^1(\mu_j)$ for all $j \in \{1,...,n\}.$ The following assertion hold:
\begin{itemize}
\item[{}]
  Assume that  $R_j:X_j \to X_j$ is  a  $\mu_j$-measure  preserving map for every $j\in \{1,...,n\}.$  If
\[c(x_1, x_2,..., x_n)=c(R_1 x_1, R_2 x_2,..., R_n x_n), \qquad \quad \forall (x_1,...,x_n) \in X_1 \times X_2 \times... \times X_n, \] then
 $(DK)$ has  a solution $(\psi_1,...,\psi_n)$  such that $\psi_j=\psi_j \circ R_j,$  $\mu_j$-a.e. on  $X_j$ with 
 \begin{eqnarray} \label{kdp}\psi_j \in {\mathcal L}^1(\mu_j) \quad \& \quad
\psi_j(x_j)=\inf \big\{ c(x_1, x_2,...,x_n) - \sum_{i=1, i \not=j}^n\psi_i (x_i);  x_i \in X_i \, \& \,i \not=j\big \},
   \end{eqnarray}
for every $j\in \{1,...,n\}$.\\
Moreover, if  $R_1, R_2,..., R_n$  are periodic -not necessary with the same period-
then $\psi_j=\psi_j \circ R_j$ on entire  $X_j.$
\end{itemize}
\end{theorem}
In case where $X_1=X_2=...=X_n,$ we have the following result.
\begin{theorem}\label{main2}
Let $X$ be a Polish space,  $\mu_1,..., \mu_n$ be n  Borel probability measures on $X,$ and
  $c: X^n  \to [0, \infty]$ be  Borel measurable and $\otimes_{i=1}^n\mu_i$-a.e.
finite. Let  $\sigma : X^n \to  X^n $  be  a  permutation defined by $\sigma(x_1,...,x_n)=(x_2,...,x_n,x_1)$. We also assume  that there exists a finite transport plan, 
 and  that  the dual problem $(DK)$ has  a solution $(\phi_1,...,\phi_n)$
such that $\phi_j \in {\mathcal L}^1(\mu)$ for all $j \in \{1,...,n\}.$ The following assertions hold:
\begin{itemize}
\item[{i.}] Suppose that   $R : X \to X$ is a   periodic map of order $n.$ If $\mu_j=(R^{n+1-j})_\# \mu_1$ for every $1 \leq j\leq n$ and   \[c(x_1, x_2,..., x_n)=c(\sigma(Rx_1, R x_2,..., Rx_n)), \qquad \quad \forall (x_1,...,x_n) \in X^n,\] then $(DK)$ has a solution $(\psi_1,...,\psi_n)$
  satisfying (\ref{kdp}) such that
    $\psi_j=\psi_1 \circ R^{j-1}$  on  $X.$

  \item[{ii.}] Suppose that $\mu_j=\mu_1$ for every $2\leq j\leq n.$ Let 
  $\{R_k\}_{(1\leq k \leq l)}$ be   a finite sequence of   $\mu_1$-measure  preserving maps on $X$ such that
 \[c(x_1, x_2,..., x_n)=c(R_k x_1, R_k x_2,..., R_k x_n), \qquad \quad \forall (x_1,...,x_n) \in X^n \, \& \,  \forall k \in \{1,...,l\}. \]
If $c$ is invariant under the permutation $\sigma$ and  $R_i \circ R_j=R_j \circ R_i $ for all $i,j \in\{1,...,l\},$ then
 $(DK)$ has  a solution $(\psi_1,...,\psi_n)$  satisfying (\ref{kdp})  such that
\[\psi_1=\psi_2=...=\psi_n, \]  and $\psi_1 \circ R_k=\psi_1,$ $\mu_1$-a.e. on $X$ for every $k \in \{1,...,l\}$.\\
 Moreover, if  $R_1, R_2,..., R_l$  are periodic -not necessary with the same period-
then $\psi_1 \circ R_k=\psi_1$ on entire  $X.$
  \end{itemize}
\end{theorem}

\begin{remark} The arguments presented in the proof of Theorems \ref{main}
 and \ref{main2} are  broad enough to capture more invariance properties  of  the dual problem. We made an effort to state
 the above theorems for a general setting. In  specific cases one may be able  to obtain more information about the dual problem.
For instance, part (ii) of Theorem \ref{main2} can be generalized to an infinite sequence of $\{R_k\}_{(1\leq k < \infty)}$ if
 solutions of the dual problem are uniformly bounded in certain  spaces.
   \end{remark}

  If the dual problem has a unique solution then  the aforementioned  Theorems will be quite useful to characterize it.  Even though the existence of an optimal transport
plan for the Monge-Kantorovich problem holds under rather general hypotheses on the cost function and marginal measures,   its uniqueness
remains an issue. In general, there are very limited cases that one can obtain  uniqueness for $(MK)$. However, the dual problem seems to have
 a better chance to admit a unique solution even when  solutions of the primal problem $(MK)$ fail to be  unique. 
  We  have the following definition.
\begin{definition}
 A Borel measure  $\mu$ on $\R^d$ is said to have a regular support if there exits a connected  open set $O$ in $\R^n$ so that $\mu(O)=1$ and
$\mu (\partial O)=0$ where $\partial O$ stands for  the boundary of $O.$
\end{definition}
 
   Here, we state the following uniqueness result for the dual problem.
\begin{proposition} \label{uniq} Let   $c: (R^d)^n \to \R$ be  differentiable and locally Lipschitz, and  $\mu_1,...,\mu_n$ be $n$  probability measures
on $\R^d$ with bounded and regular supports.    If each  $\mu_i$ is absolutely continuous with
respect to the  $d$-dimensional Lebesgue measure, ${\mathcal L}^d,$  and has a positive density on its support  with respect to ${\mathcal L} ^d$, then $(DK)$ admits
a unique solution (up to the  addition of constants summing to $0$ to each potential).
\end{proposition}
Similar uniqueness results for the problem  $(DK)$ can be found in \cite{G-O-Z, C-M-N, Lo}. However, for the convenience of the reader we shall  provide a proof in this paper.
Finally, we shall show that an optimal plan can be chosen so that it inherits the  invariance properties of   the cost function and the marginal measures.
\begin{theorem} \label{MK} Let  $(X_1, \mu_1),...,(X_n, \mu_n)$ be $n$  Polish probability   spaces,
and  $c: X_1 \times X_2 \times... \times X_n  \to [0, \infty]$ be  a  lower  semi-continuous Borel measurable function.
 The following statements hold.

\begin{itemize}
\item[{i.}] For each $j\in \{1,...,n\},$ let $R_j:X_j \to X_j$ be a  $\mu_j$-measure  preserving map. If
\[c(x_1, x_2,..., x_n)=c(R_1 x_1, R_2 x_2,..., R_n x_n), \qquad \quad \forall (x_1,...,x_n) \in X_1 \times X_2 \times... \times X_n, \] then
$(MK)$ has a solution $\bar \gamma$ such that $ \bar \gamma=(R_1,..., R_n)\#\bar \gamma.$
\item[{ii.}] Let $X_1=...=X_n=X$ and $\mu_1=...=\mu_n=\mu.$ Let $\mathcal {G}$ be a set of
   $\mu$-measure  preserving maps on $X$ such that
 \begin{equation*}c(x_1, x_2,..., x_n)=c(U x_1, U x_2,..., Ux_n), \qquad \quad \forall (x_1,...,x_n) \in X^n \, \& \,  \forall U  \in \mathcal {G}, \end{equation*}
and,  that  $U \circ R=R \circ  U $ for all $ U, R  \in \mathcal {G}.$  If $\mathcal {G}$ is a countable set  then
$(MK)$ has a solution $\bar \gamma$ that is invariant under  $\mathcal {G},$  i.e. $ \bar \gamma=(U,...,U)\#\bar \gamma$ for all $ U  \in \mathcal {G}.$\\
  \end{itemize}
\end{theorem}
The following result is an immediate  consequence of the above theorem.
\begin{coro}\label{mkc} Under the assumptions of part (ii) in Theorem \ref{MK}, if the cost function $c$ is invariant by a
 permutation $\sigma : X^n \to X^n,$ i.e. $c \circ \sigma=c$ then
$(MK)$ has a solution $\bar \gamma$ such that  $\bar \gamma=\sigma \# \bar \gamma,$ and  $ \bar \gamma=(U,..., U)\#\bar \gamma$ for all $U  \in \mathcal {G}.$\\
Furthermore, If  $X$ is    a metric space  and $\mathcal {G}$ equipped with the point-wise topology on $X$  is separable then
 the  countability  assumption on  $\mathcal {G}$ is not required.
 \end{coro}
\begin{remark} By drawing a comparison between  Theorems  \ref{main2}  and \ref{MK}, it is evident that the primal problem $(MK)$ is more flexible than  its dual $(DK)$ when it
comes
to capturing invariance  properties of the cost function and marginal measures. The reason falls under the fact that transport plans live in the set $\Pi(\mu_1,...,\mu_n)$
 which is  pre-compact for the weak topology. It allows one to analyze the limit points of optimal plans with certain properties. \end{remark}
The literature on the theory and applications of optimal mass transportation  is too vast to provide an exhaustive bibliography here,  we refer to the books
of Villani \cite{V} and Rachev and R\"uchendrof \cite{RR}
and the references therein.\\

The next section is devoted to the proof of Theorems \ref{main} and \ref{main2}. In section 3, by recalling the measure isomorphism theorem we
 address the uniqueness issue for the dual problem. Section 4 is concerned with the invariance properties of $(MK)$ and the proof of Theorem \ref{MK}.
Applications to optimal mass transportation problems arising in volume maximizing and also   semi-classical Hohenberg-Kohn functionlas are discussed in section 5.

\section{ Invariance properties of $(DK)$}

This section is devoted to the proof of Theorems \ref{main} and \ref{main2}. We first recall the celebrated Ergodic theorem, due to George David Birkhoff (1931).
\begin{theorem} \label{ergodic}{\bf(Birkhoff's Ergodic Theorem)} Let $T: X \to X$ be a measure-preserving transformation on a measure space $(X, \Sigma,\mu)$
 and suppose $f$ is a $\mu$-integrable function, i.e. $f \in {\mathcal L}^1(\mu).$ Then the average
\[A_m(f(x))= \frac{1}{m} \sum_{i=0}^{m-1} f(T^ix),\]
converge a.e. on $X$ (as $m \to \infty$) to an integrable function $\hat f.$
Furthermore, $\hat f $ is $T$-invariant, i.e.
   $\hat f \circ T= \hat f \,$
holds almost everywhere, and if $\mu (X)$ is finite, then the normalization is the same,
    \[\int_X \hat f\, d\mu = \int_X f\, d\mu.\]
\end{theorem}
The above theorem is an essential tool in Ergodic theory to study  the behavior of  a dynamical system  with an invariant measure. We shall see that it is also
surprisingly relevant and essential in the theory of optimal mass transportation. \\
\textbf{Proof of Theorem \ref{main}.}
For each $j\in \{1,...,n\},$  define the average of $\phi_j$  under the the transformation $R_j$ as in Theorem \ref{ergodic},
\[A_m (\phi_j(x))=\frac{\sum_{k=0}^{m-1}\phi_j (R^k_j x)}{m}.\]
It follows from Theorem \ref{ergodic} that $A_m (\phi_j)$ converges $\mu_j$-a.e.  to an integrable function. Define
 \begin{eqnarray*}
\Phi_j(x)=\left\{
\begin{array}{ll}
\lim_{ \substack{ m \to \infty}}A_m \big (\phi_j(x) \big), & \text{ provided this limit exists, } \\
-\infty, & \text{ otherwise. }\\
\end{array}
\right.
\end{eqnarray*} 
It follows that  $\Phi_j$ must satisfy 
\begin{eqnarray}\label{erg}\Phi_j \circ R_j= \Phi_j,\,  \quad  \, \mu_j- a.e. \, \qquad \&  \qquad \int_{X_j}\Phi_j \, d \mu_j= \int_{X_j}\phi_j \, d \mu_j.\end{eqnarray}
 Denote by $\Phi_j^c$ the c-conjugate of $\Phi_j$ defined by
\[\Phi_j^c(x_j)=\inf \Big \{ c(x_1,...,x_n) - \sum_{i=1, i \not=j}^n\Phi_i (x_i); \,
x_i \in X_i \, \& \, i \not=j\Big \}.\]
{\it Claim.} For each $j \in \{1,...,n\},$ $\Phi_j^c$ satisfies the following properties:\\
a.   $\Phi_j\leq \Phi^c_j$ on $X_j.$\\
b. $\Phi_j^c \in {\mathcal L}^1(\mu_i).$\\
c.  $\Phi_j^c =\Phi_j,$  $\mu_j$-a.e. on  $ X_j.$\\

 Let us first prove the claim.  Fix  $(x_1,...,x_n) \in \tilde X= X_1 \times X_2\times... \times X_n$ such that  $\lim_{ \substack{ m \to \infty}}A_m \big (\phi_i(x_i) \big)$ exist
 for every  $i.$
 We have,
\begin{eqnarray*}
 c(x_1,...,x_n)-\sum_{i=1, i \not=j}^n\Phi_i (x_i)&=& c(x_1,...,x_n)-\lim_{m \to \infty }\sum_{i=1, i \not=j}^n\frac{\sum_{k=0}^{m}\phi_i(R^k_i x_i)}{m}\\
&=&\lim_{m \to \infty }\sum_{k=0}^{m}\frac{c(R_1^k x_1,...,R_n^k x_n)-\sum_{i=1, i \not=j}^n\phi_i(R^k_i x_i)}{m}\\
&\geq&\lim_{m \to \infty }\frac{\sum_{k=0}^{m} \phi_j(R^k_i x_j)}{m}=\Phi_j(x_j),\\
\end{eqnarray*}
 and therefore \[c(x_1,...,x_n)\geq \sum_{i=1}^n\Phi_i (x_i).\] Note also that the latter inequality still holds if   $\lim_{ \substack{ m \to \infty}}A_m \big (\phi_i(x_i) \big)$ 
does not exist for some $i,$ (since  the right  hand side will be $-\infty$).
It then follows that   $\Phi_j\leq \Phi_j^c$.  To prove part (b), note that
\[\Phi_j(x_j)   \leq \Phi_j^c(x_j)\leq c(x_1,...,x_n) - \sum_{i=1, i \not=j}^n\Phi_i (x_i), \qquad \forall (x_1,...,x_n) \in \tilde X. \]
and therefore
\[  |\Phi_j^c(x_j)| \leq c(x_1,...,x_n) + \sum_{i=1}^n|\Phi_i (x_i)|, \qquad \forall (x_1,...,x_n) \in \tilde X. \]
By Birkhoff's Ergodic theorem  each $\Phi_i$ is integrable and by the assumption  there exists a finite transport plan  from which   the integrability of $ \Phi_j^c$ follows.
To prove part (c), we first note that
 \begin{eqnarray*} \Phi_j^c(x_j)+ \sum_{i=1, i \not=j}^n\Phi_i (x_i) \leq   c(x_1, x_2,...,x_n), \qquad \quad \forall (x_1,...,x_n) \in \tilde X.
   \end{eqnarray*}
It  follows from part (a) of the claim and (\ref{erg}) that
 \begin{eqnarray*} \int_{X_j}\Phi_j^c(x_j)\, d\mu_j+ \sum_{i=1, i \not=j}^n\int_{X_i}\Phi_i (x_i) \, d\mu_i\geq  \sum_{i=1}^n \int_{X_i} \Phi_i (x_i) \, d\mu_i =
  \sum_{i=1}^n \int_{X_i} \phi_i (x_i) \, d\mu_i.
   \end{eqnarray*}
The optimality of $(\phi_1,...,\phi_n)$ implies that the inequality in the latter expression is  in fact an equality.
Therefore, $\int_{X_j}\Phi_j^c(x_j)\, d\mu_j= \int_{X_j}\Phi_j(x_j)\, d\mu_j$ from which together with the fact that $\Phi_j^c \leq \Phi_j$
 we obtain $\Phi_j^c = \Phi_j,$ $\mu_j$-a.e. on $X_j.$  This also implies that    $\Phi_j^c \circ R_j = \Phi_j^c,$ $\mu_j$-a.e.  on  $X_j$ as the same property holds by $\Phi_j$.\\
 Note that $(\Phi_1,...,\Phi_n)$ satisfies the required invariance property.  It also satisfies the $c$-concavity condition (\ref{kdp}) up to a null set as 
$\Phi_j^c = \Phi_j,$ $\mu_j$-a.e. on $X_j.$  Following an idea in \cite{G-S, R}, we shall now show that $(\Phi_1,...,\Phi_n)$ can be replaced by an $n$-tuple $(\psi_1,...,\psi_n)$ satisfying (\ref{kdp}) on
 the whole $\tilde X.$  

 It  follows from the definition of $\Phi_j^c$ that,
\[\sum_{i=1}^n\Phi_i^c(x_i)+ (n-1) \sum_{i=1}^n\Phi_i(x_i) \leq n c(x_1,x_2,...,x_n), \qquad\qquad  \forall(x_1, x_2, ..., x_n) \in \tilde X. \]
Setting \[V_i(x_i)= \frac{\Phi_i^C(x_i)+(n-1)\Phi_i(x_i)}{n}, \qquad x_i \in X_i,\] the latter  inequality yields that
\[\sum_{i=1}^n V_i(x_i)\leq c(x_1,...,x_n), \qquad\qquad  \forall(x_1, x_2, ..., x_n) \in \tilde X.\]
 It also  follows from the above claim that $\Phi_j \leq V_j \leq \Phi_j^c $ on $X_j$ and $V_j \circ R_j=V_j,$  $\mu_j$-a.e. on $X_j.$
Let now
\begin{eqnarray*} \psi_1(x)=\sup\Big \{w(x);  \, \Phi_1  \leq w \leq \Phi_1^c \,  \& \,
w(x_1)+ \sum_{i=2}^n V_i(x_i)\leq c(x_1,...,x_n),
 \, \forall (x_1,x_2,...,x_n)  \in  \tilde X\Big \}.
\end{eqnarray*}
By induction
 for each $1< k < n$ define,
\begin{eqnarray*}\psi_k(x)=\sup \Big \{w(x); \,  \Phi_k  \leq w \leq \Phi_k^c \, \& \,
 \sum_{i=1}^{k-1}\psi_i(x_i)+ w(x_k)+  \sum_{i=k+1}^n V_i(x_i)\leq c(x_1,...,x_n), \, \text{on} \, \tilde X
\Big \},\end{eqnarray*}
and finally we  define
\begin{eqnarray*}\psi_n(x)=\sup\Big \{w(x);  \, \Phi_n  \leq w \leq \Phi_n^c \, \& \,
\sum_{i=1}^{n-1}\psi_i(x_i)+ w(x_n)\leq c(x_1,...,x_n), \, \text{on} \, \tilde X
\Big \}.\end{eqnarray*}
Note that each $\psi_j$ is measurable with respect to the $\mu_j-$completion of the Borel $\sigma-$algebra on $X_j$. Indeed, the measurability  follows from the inequalities  $\Phi_j \leq \psi_j\leq \Phi_j^c$  and the fact that $\Phi_j=\Phi_j^c,$ $\mu_j-$a.e.\\

 It follows from the definition of $\psi_j$ that,
\begin{equation} \label{ine4}
\Phi_j^c \geq \psi_j \geq V_j  \geq \Phi_j,
\end{equation}
and,
\begin{equation} \label{ine3}
\sum_{i=1}^{n}\psi_i(x_i)\leq c(x_1, x_2,...,x_n), \qquad \forall   (x_1,x_2,...,x_n)  \in  X.
\end{equation}
Consider now the functions,
\begin{eqnarray} \label{ine1} \bar \psi_j(x_j)=\inf \Big \{c(x_1, x_2,...,x_n)-\sum_{i=1, i \not=j}^n \psi_i(x_i); \, x_i \in X_i \, \& \, i \not =j \Big\}.\end{eqnarray}
It follows from (\ref{ine3}) that  $\bar \psi_j \geq \psi_j.$ It also follows from (\ref{ine4}) that
\begin{eqnarray*}
\bar \psi_j(x_j) &\leq& \inf \Big \{c(x_1,...,x_n)-\sum_{i=1, i \not=j}^n \Phi_i(x_i); \, x_i \in X_i \, \& \, i \not =j \Big\}=\Phi_j^c(x_j),
\end{eqnarray*}
and therefore,
\begin{eqnarray} \label{ine2}  \Phi_j^c \geq \bar \psi_j \geq \psi_j \geq \Phi_j.\end{eqnarray}
It then   follows from (\ref{ine1}), (\ref{ine2})  and  the maximality of $\psi_j$ that $\bar \psi_j=\psi_j$
for all $j\in \{1,...,n\}.$ Therefore, one obtains
 \begin{eqnarray*}  \psi_j(x_j)=\inf \Big \{c(x_1, x_2,...,x_n)-\sum_{i=1, i \not=j}^n \psi_i(x_i); \, x_i \in X_i \, \& \, i \not =j \Big\}.
   \end{eqnarray*}
 To complete the proof we will show that $(\psi_1,\psi_2,...,\psi_n)$ is a solution of $(DK).$  The inequality  relation  in (\ref{ine2}) together with
 the integral equality given in (\ref{erg}) yield that
 \begin{eqnarray*}
\sum_{i=1}^n\int_{X_i} \psi_i(x_i) \, d \mu_i \geq \sum_{i=1}^n\int_{X_i} \Phi_i (x_i) \, d \mu_i
= \sum_{i=1}^n\int_{X_i} \phi_i(x_i) \, d \mu_i,
\end{eqnarray*}
and therefore, the optimality  of $(\phi_1, \phi_2,...,\phi_n)$  implies that
\[\sum_{i=1}^n\int_{X_i} \psi_i(x_i) \, d \mu_i=\sum_{i=1}^n\int_{X_i} \phi_i(x_i) \, d \mu_i.\]
This completes the proof of part (i) for possibly non-periodic maps $R_1,..., R_n$. Now assume that $R_1,..., R_n$ are periodic.
There exist positive integers $m_1,..., m_n$ such that
$R_j^{m_j}$ is the identity map on $X_j$  for each $j \in \{1,...,n\}.$  Define  $\Phi_j$ to be
\[\Phi_j (x)=\frac{\sum_{k=0}^{m_j-1}\phi_j (R^k_j x)}{m_j}, \qquad \forall x \in X_j.\]
Note that $\Phi_j=\Phi_j \circ R_j$ on entire $X_j.$  One can now  deduce that $\Phi^c_j=\Phi^c_j \circ R_j$ on entire $X_j.$ In fact,
\begin{eqnarray*}\Phi_j^c(x_j)&=&\inf \{ c(x_1, x_2,...,x_n) - \sum_{i=1, i \not=j}^n\Phi_i (x_i); \,
x_i \in X_i \, \& \, i \not=j\}\\
&=&\inf \{ c(R_1x_1, R_2x_2,...,R_nx_n) - \sum_{i=1, i \not=j}^n\Phi_i (R_ix_i); \,
x_i \in X_i \, \& \, i \not=j\}\\
&=& \Phi_j^c(R_j x_j), \qquad (R_j \text{ is surjective as } R_j^{m_j}=Id).
\end{eqnarray*}
The rest of the proof goes in the same lines as in the previous case. \hfill $\square$\\

For the case where $X_1=...=X_n=X$ one can hope for more invariance properties as stated in Theorem \ref{main2}. We shall now proceed with the proof of this Theorem.\\

\textbf{Proof of Theorem \ref{main2}.}
Define, \[\Phi (x)=\frac{\sum_{k=1}^{n}\phi_k (R^{n-k} x)}{n},\] and
$\Phi_i(x)=\Phi(R^{i} x)$  for $i \in \{1,...,n\}$.  It can be easily deduced that
\begin{equation}\label{late}
\int_X \Phi_i(x) \, d\mu_i=\frac{1}{n} \sum_{k=1}^n \int_X \phi_k(x) \, d\mu_k
 \end{equation}
for each $i \in \{1,...,n\}$. Denote by $\Phi_1^c$ the c-conjugate of $\Phi_1$ defined by
\[\Phi_1^c(x_1)=\inf \big \{ c(x_1, x_2,...,x_n) - \sum_{i=2}^n\Phi_i (x_i);   x_2,...,x_n  \in X \big\}.\]
Setting \[V(x)= \frac{\Phi_1^c(x)+(n-1)\Phi_1 (x)}{n},\]it follows from the above expression that
\[\sum_{i=1}^n V(R^{i-1}x_i)\leq c(x_1,...,x_n), \qquad \forall(x_1, x_2, ..., x_n) \in  X^n.\]
We also have that $\Phi_1\leq \Phi_1^c.$ In fact,
\begin{eqnarray*}
 c(x_1,...,x_n)-\sum_{i=2}^n\Phi_i (x_i)&=&c(x_1,...,x_n)-\sum_{i=2}^n\Phi(R^{i} x_i) \\
&=&c(x_1,...,x_n)-\frac{\sum_{i=2}^n\sum_{k=1}^{n} \phi_k(R^{n-k} R^ix_i)}{n}\\
&=&c(x_1,...,x_n)-\frac{\sum_{i=2}^n\sum_{j=1}^{n} \phi_{n+i-j}(R^jx_i)}{n}, \qquad  \big (\phi_{n+k}=\phi_k, \, \, \forall k\in \mathbb{N}\big ),\\
&=&\sum_{j=1}^{n}\frac{c\big( \sigma^{j}(R^{j} x_1,...,R^{j} x_n) \big)-\sum_{i=2}^n\phi_{n+i-j}(R^jx_i)}{n}\\
&\geq&\frac{\sum_{j=1}^{n} \phi_{n+1-j}(R^jx_1)}{n}=\frac{\sum_{k=1}^{n} \phi_{k}(R^{n+1-k}x_1)}{n}=\Phi_1(x_1)\\
\end{eqnarray*}
from which one has $\Phi_1\leq \Phi_1^c$ . It follows that $\Phi_1^c(x) \geq V(x) \geq \Phi_1(x) $ for every $x \in X.$
Let now
\begin{eqnarray*} \psi(x)=\sup\Big \{w(x);  \quad  \Phi_1  \leq w \leq \Phi_1^c \,  \& \,
 \sum_{i=1}^n w(R^{i-1}x_i)\leq c(x_1,...,x_n),
 \,\quad  \forall (x_1,x_2,...,x_n)  \in  X^n\Big \}.
\end{eqnarray*}
Consider now the function,
\begin{eqnarray*}  \bar \psi(x_1)=\sup \Big \{c(x_1, x_2,...,x_n)-\sum_{i=2}^n \psi(R^{i-1}x_i); \, x_i \in X \, \& \, i \not =j \Big\}.\end{eqnarray*}
It follows that $\bar \psi \geq \psi.$
Set  $w(x)=\frac{\bar \psi (x) +(n-1) \psi(x)}{n}$ and note that $w \geq \psi.$   One can easily check that
 \[\Phi_1^c \geq w \geq \Phi_1 \quad \& \quad \sum_{i=1}^n w(R^{i-1}x_i)\leq c(x_1,...,x_n)\quad \text{ on }X^n.\]
By the  maximality  of $\psi$ we have that $\psi=w$ and therefore
  $\bar \psi = \psi.$  We now show that $(\psi, \psi \circ R,..., \psi\circ R^{n-1})$ is a solution of $(DK).$ Since $\Phi_1 \leq \psi$ we have that 
$\Phi_i \leq \psi \circ R^{i-1}$ from which we obtain
\begin{equation}\label{late2}
\sum_{i=1}^n \int_X\Phi_i(x) \, d\mu_i \leq \sum_{i=1}^n \int_X\psi(R^{i-1}x) \, d\mu_i.
\end{equation}
On the other hand by (\ref{late}) we have that $\sum_{i=1}^n \int_X\Phi_i(x) \, d\mu_i =\sum_{i=1}^n \int_X\phi_i(x) \, d\mu_i$ and therefore it follows from (\ref{late2}) that  
\[\sum_{i=1}^n \int_X\phi_i(x) \, d\mu_i \leq \sum_{i=1}^n \int_X\psi(R^{i-1}x) \, d\mu_i.\]
 This together with the  optimality of $(\phi_1,...,\phi_n)$ 
completes the proof of part (i).\\

  Proof of part (ii). By assuming $R$ to be the identity map in part (i), it follows  that a solution $(\phi_1,...,\phi_n)$ of $(DK)$ can be chosen in such a
 way that $\phi_j=\phi_1$ for all $j \in \{1,...,n\}.$   Define the average of $\phi_1$  under the the transformation $R_1$ as in Theorem \ref{ergodic},
\[A_m (\phi_1(x))=\frac{\sum_{k=0}^{m-1}\phi_1 (R^k_1 x)}{m}.\]
It follows from Theorem \ref{ergodic} that $A_m (\phi_j)$ converges $\mu$-a.e.    to an integrable function $\Phi_1$ satisfying
\begin{eqnarray*}\Phi_1 \circ R_1= \Phi_1,\,  \quad  \, \mu- a.e. \, \qquad \&  \qquad \int_{X}\Phi_1 \, d \mu_1= \int_{X}\phi_1 \, d \mu,\end{eqnarray*}
 where $\Phi_1$ is defined to be $-\infty$ at points where $A_m (\phi_j)$ does not converge.  Note also that
\[\sum_{i=1}^n A_m (\phi_1(x_i))=\frac{\sum_{k=0}^{m-1}\sum_{i=1}^n\phi_1 (R^k_1 x_i)}{m} \leq \frac{\sum_{k=0}^{m-1}c(R_1^k x_1,...,R_1^kx_i)}{m}=c(x_1,...,x_n),\]
from which we obtain
\begin{equation}\label{org}
\sum_{i=1}^n \Phi_1(x_i) \leq c(x_1,...,x_n), \qquad \qquad \forall (x_1,...,x_n)\in X^n.
\end{equation}
We now define the average of $\Phi_1$  under the the transformation $R_2$,
\begin{equation}\label{org1} A_m (\Phi_1(x))=\frac{\sum_{k=0}^{m-1}\Phi_1 (R^k_2 x)}{m}.\end{equation}
 By making use of Theorem \ref{ergodic} again, it follows that $A_m (\Phi_1)$ converges $\mu$-a.e.    to an integrable function $\Phi_2$ satisfying
\begin{eqnarray*}\Phi_2 \circ R_2= \Phi_2,\,  \quad  \, \mu- a.e. \, \qquad \&  \qquad \int_{X}\Phi_2 \, d \mu= \int_{X}\Phi_1 \, d \mu,\end{eqnarray*}
and $\Phi_2$ is defined to be $-\infty$  at points where  $A_m (\Phi_1)$ does not converge. 
 Note that $R_1 \circ R_2=R_2 \circ R_1$ and $\Phi_1$ is invariant under $R_1,$  $\mu$-a.e. on $X.$ Thus, $A_m (\Phi_1(R_1x))=A_m (\Phi_1(x))$ for $\mu$-a.e.  $x\in X.$
 Since $A_m (\Phi_1)$ converges $\mu$-a.e.  to $\Phi_2$ it follows that $\Phi_2 \circ R_1= \Phi_2,$  $\mu$-a.e. on $X.$ It also follows from (\ref{org}) and (\ref{org1}) that
 \begin{equation*}
\sum_{i=1}^n \Phi_2(x_i) \leq c(x_1,...,x_n), \qquad \qquad \forall (x_1,...,x_n)\in X^n.
\end{equation*}
By repeating the above argument we obtain, by the {\it l}-th step, a function  $\Phi_l$ enjoying the following properties:
1. $\int_{X}\Phi_l \, d \mu=...=\int_{X}\Phi_1 \, d \mu=\int_{X}\phi_1 \, d \mu.$\\
2. $ \Phi_l \circ R_k=\Phi_l,$ $\mu$-a.e.  on $X$  for each $k \in \{1,...,l\}.$\\
3. $\sum_{i=1}^n \Phi_l(x_i) \leq c(x_1,...,x_n),  \quad \forall (x_1,...,x_n)\in X^n.$\\

Denote by $\Phi_l^c$ the c-conjugate of $\Phi_l$ defined by
\[\Phi_l^c(x_1)=\inf \big \{ c(x_1, x_2,...,x_n) - \sum_{i=2}^n\Phi_l (x_i);   x_2,...,x_n  \in X \big\}.\]
Setting $V(x)= \frac{\Phi_l^c(x)+(n-1)\Phi_l (x)}{n}$, it follows from the above expression that
\[\sum_{i=1}^n V(x_i)\leq c(x_1,...,x_n), \qquad \forall(x_1, x_2, ..., x_n) \in  X^n.\]
We  have that $\Phi_l^c\geq \Phi_l$  as $\sum_{i=1}^n \Phi_l(x_i) \leq c(x_1,...,x_n)$ on $X^n.$ We also have that $\Phi_l^c(x)= \Phi_l(x)$ for $\mu$-a.e. $x \in X$ (as in the part (c)
 of the claim in the proof of Theorem \ref{main}). It follows that $\Phi_l^c(x) \geq V(x) \geq \Phi_l(x) $ for each $x \in X.$
Let now
\begin{eqnarray*} \psi(x)=\sup\Big \{w(x);  \quad  \Phi_l \leq w \leq \Phi_l^c \,  \& \,
 \sum_{i=1}^n w(x_i)\leq c(x_1,...,x_n),
 \,\quad  \forall (x_1,x_2,...,x_n)  \in  X^n\Big \}.
\end{eqnarray*}
Similar to part (i), one has
\begin{eqnarray*}   \psi(x_1)=\inf \Big \{c(x_1, x_2,...,x_n)-\sum_{i=2}^n \psi(x_i); \, x_i \in X  \Big\}.\end{eqnarray*}
Note that   $\Phi_l^c  \geq \psi \geq \Phi_l,$ and  $\Phi_l$ and $\Phi_l^c$ are $R_k$ invariant $\mu$-a.e. on $X.$ It follows that
$\psi\circ R_k=\psi,$ $\mu$-a.e. on $X$ for all $k \in \{1,...,l\}.$
This completes the proof for possibly non-periodic maps $R_1,...,R_l.$ A similar argument as in the proof of Theorem \ref{main}
shows that if $R_1,...,R_l$ are periodic then $\psi\circ R_k=\psi$ on entire $X.$
\hfill $\square$\\

\section{Uniqueness}
In this section we  address the uniqueness issues for both $(MK)$ and $(DK)$.  Under rather general  conditions on the cost and marginals,
existence of  both $(MK)$ and $(DK)$ are warranted. For $(MK)$ with $n=2$,  the well known twist condition i.e.,
\[x_2\to  D_{x_1}c(x_1, x_2) \qquad \text{is injective for fixed} \, x_1, \]
ensures the uniqueness and
Monge structure of the optimal map \cite{G-Mc, L}. For larger $n$, the uniqueness  question
is  still largely open. Examples of special cost functions for which the optimal
measure has this structure are known \cite{C, G-S, H, H-K, L, P}  as well as several
examples for which uniqueness and Monge solutions fail \cite{C-N, P, G-M0}.

 We now proceed with the proof of Proposition \ref{uniq}.
To do this, we first recall some preliminary notations and results as in the theory of measure isomorphisms. Let  $X$ be  a topological
space and  denote by $B(X)$ the set of all Borel subsets of $X$.
\begin{definition} Assume that $\mu$ is a Borel measure on a topological
space $X$ and $\nu$ is a Borel measure on a topological space $Y$. We
say that $(X,B(X), \mu)$ is isomorphic to $(Y,B(Y ), \nu)$ if there exists a
one-to-one map $T$ of $X$ onto $Y$ such that for all $A \in B(X)$ we have
$T(A) \in B(Y)$ and $\mu[A] = \nu[T(A)],$ and for all $B \in B(Y)$ we have
$T^{-1}(B) \in B(X)$ and $\mu[T^{-1}(B)] = \nu[B]$.
\end{definition}
The following is the standard measure isomorphism theorem.
\begin{theorem}\label{iso}
Let $\mu$ be a finite Borel measure on a complete separable
metric space $X$. Assume that $\mu$  is non-atomic  and $\mu[X] = 1$. Then
$(X, B(X),  \mu)$ is isomorphic to $([0, 1], B([0, 1]), \lambda_1)$ where $\lambda_1$ stands for
the one-dimensional Lebesgue measure on $[0, 1].$
\end{theorem}

\textbf{Proof of Proposition  \ref{uniq}} Since each  $\mu_i$ has a regular support there exist a connected open set $X_i$ with  full measure such that that
 $\mu_i(\partial X_i)=0.$ Let  $\tilde X= X_1 \times...\times X_n,$ and assume that
 $\gamma$ is a solution of $(MK).$ Since each $\mu_i$ is non-atomic we have that $\gamma$ is a non-atomic measure on $B(\tilde X)$. It follows from Theorem
\ref{iso} that $(\tilde X, B(\tilde X), \gamma)$ is isomorphic to $(X_1, B(X_1), \mu_1).$ Thus, there exists an isomorphism $T=(T_1,...,T_n)$ from $(X_1, B(X_1), \mu_1) $ onto
 $(\tilde X, B(\tilde X), \gamma)$.  Note that each $T_i: X_1 \to X_i$ is onto and  pushes $\mu_i$ forward to $\mu_1,$ i.e.
\begin{eqnarray}\label{diff}\mu_1(T_i^{-1}A)=\mu_i(A), \qquad \forall  A \in B(X_i)\, \, \&\, \, \forall i \in \{1,...,n\}. \end{eqnarray}
  Let $(\phi_1,...,\phi_n)$ be a solution of $(MK)$ satisfying
\begin{eqnarray} \label{kdp000}\phi_j(x_j)=\inf \big\{ c(x_1, x_2,...,x_n) - \sum_{i=1, i \not=j}^n\phi_i (x_i);  x_i \in X_i \, \& \,i \not=j\big \}.
   \end{eqnarray}
   Since $c$ is locally Lipschitz and $X_1,...,X_n$ are bounded, each $\phi_i$  is locally Lipschitz (Lemma C. 1 \cite{G-Mc}). Suppose
that  $A_i \subset \R^d$ is  the set of non-differentiable points of $\phi_i$ on $X_i.$  It follows from Rademacher's theorem that ${\mathcal L}^d( A_i)=0,$ and since $\mu_i$
is absolutely continuous with respect to ${\mathcal L}^d$, one has $\mu_i(A_i)=0.$  It now follows from (\ref{diff}) that $\mu_1(\cup_{i=1}^nT_i^{-1}A_i)=0.$
   For every  $x \in X_1\setminus \cup_{i=1}^nT_i^{-1}A_i,$ each $\phi_i$ is differentiable at $T_ix.$\\
 On the other hand, it follows from (\ref{diff}) and the duality between $(MK)$ and $(DK)$  that
   \begin{eqnarray*}
   \int_{X_1} c(T_1x,...,T_nx) \, d \mu_1&=&\int_{\tilde X}c(x_1,...,x_n) \, d \gamma\\
   &=&\sum_{i=1}^n\int_{X_i} \phi_i(x) \, d \mu_i\\
   &=&\sum_{i=1}^n\int_{X_i} \phi_i(T_ix) \, d \mu_1.
   \end{eqnarray*}
Since $c(x_1,...,x_n)\geq \sum_{i=1}^n \phi_i(x_i)$, there exists  a  measurable subset $A_0$ of $X_1$ with $\mu_1(A_0)=0$ such that

   \[c(T_1x,...,T_nx)=\sum_{i=1}^n \phi_i(T_ix), \qquad \, \forall  x \in X_1\setminus A_0.\]
  Set $A= \big( \cup_{i=1}^nT_i^{-1}A_i\big ) \cup A_0 $ and note that $\mu_1(A)=0.$  For each $x \in X_1\setminus A,$ the differentiability of $c$ together with the fact that $c(x_1,...,x_n)\geq \sum_{i=1}^n \phi_i(x_i)$,
yield  that
   \[\frac{\partial c(T_1x,...,T_nx)}{\partial x_i}=\nabla \phi_i(T_ix), \qquad  \, \forall x  \in X_1\setminus A.\]
If $(DK)$ has another solution $(\psi_1,...,\psi_n)$ satisfying (\ref{kdp000}), the above argument shows that
 \[\frac{\partial c(T_1x,...,T_nx)}{\partial x_i}=\nabla \psi_i(T_ix), \qquad  \, \forall x  \in X_1\setminus B,\]
for some measurable set $B$ with  $\mu_1(B)=0.$ Therefore,
\begin{equation}\label{abb}\nabla \phi_i(T_ix)=\nabla \psi_i(T_i x) \qquad \forall x \in X_1\setminus (A\cup B).\end{equation}
We show that $\nabla \phi_i(z)=\nabla \psi_i(z)$ for ${\mathcal{L}}^d$-a.e. $z \in X_i.$
Let $\lambda_i$ be the measure ${\mathcal L}^d$  restricted to $X_i$. Since $T$ is a measurable isomorphism  we have that $T\big (X_1\setminus (A\cup B)\big)$ is measurable and \[\gamma \Big(T\big (X_1\setminus (A\cup B)\big)\Big)=1.\]
Let $\Lambda_i$ be the projection of $T\big (X_1\setminus (A\cup B)\big)$ on the i-th variable.
Thus,  $\Lambda_i$ is universally measurable and there exists Borel sets $K_i, O_i$ with $K_i \subseteq \Lambda_i \subseteq O_i$ such that $\mu_i(O_i \setminus K_i)=0.$ Since
\[T\big (X_1\setminus (A\cup B)\big)\subseteq X_1\times...O_i\times...\times X_n,\]
we obtain that $\mu_i(O_i)=1.$ Thus, $\mu_i(K_i)=1.$ Note also that since $K_i \subseteq \Lambda_i$ it follows from (\ref{abb}) that 

\begin{equation}\label{abb1}\nabla \phi_i(z)=\nabla \psi_i(z) \qquad \forall z \in K_i.\end{equation}
  Since by assumption $\frac{d \mu_i}{ d \lambda_i}$ is a  positive function and $\mu_i(K_i)=\mu_i(X_i)=1,$ we obtain that  $\lambda_i(K_i)=\lambda_i(X_i).$ It now follows from (\ref{abb1}) that 
 \begin{equation*}\nabla \phi_i(z)=\nabla \psi_i(z) \qquad {\mathcal{L}}^d-a.e. \, z \in X_i,\end{equation*}
 from which the desired result follows.  \hfill $\square$\\

\section{Invariance properties of $(MK)$}

We first recall the following result from  (\cite{V}, Lemma 4.3).
\begin{lemma}\label{vilan}
 Let $X_1,...,X_n$  be Polish spaces and $c: X_1 \times...\times X_n \to [0,\infty]$ be lower  semi-continuous. Assume that $\gamma_k$ is a sequence in  $\Pi(\mu_1,...,\mu_n)$
converging weakly to a transport plan $\gamma.$ Then
\[\int c \, d \gamma \leq \liminf_{k \to \infty}\int c \, d \gamma_k.  \]
\end{lemma}

The above result does not imply that the optimal cost is finite. In fact, all transport plans may lead to an infinite cost i.e. $\int c \, d\gamma=\infty$ for all
 $\gamma \in \Pi(\mu_1,...,\mu_n).$  We are now ready to prove  Theorem \ref{MK}.\\

\textbf{Proof of Theorem \ref{MK}.}  Let us first  assume that there exists a finite transport plan.  Let $\tilde X= X_1 \times...\times X_n,$ and define
 $R: \tilde X  \to \tilde X $  by \[R(x_1,...,x_n)=(R_1 x_1,...,R_n x_n).\]

For each nonnegative integer $k,$ set $R^k=(R^k_1,...,R^k_n)$  with the convention that $R_i^0$ to be the identity map  on $X_i$ for each $i \in \{1,...,n\}.$
 By Lemma \ref{vilan} and the fact that there exists a finite transport plan, the existence of  an optimal plan $\gamma$ with a finite cost is ensured.  Define
\[\gamma_k=\frac{\sum_{i=0}^k R^i\# \gamma }{k+1}.\]
Note first that $\gamma_k \in \Pi (\mu_1,...,\mu_n).$  Indeed, if $f$ is a bounded continuous function on $X_j$ then
\begin{eqnarray*}
\int_{\tilde X} f(x_j) \, d \gamma_k=\frac{\sum_{i=0}^k \int_{\tilde X} f(x_j) \, d( R^i\# \gamma )}{k+1}&=&\frac{\sum_{i=0}^k \int_{\tilde X} f(R_j^i x_j) \, d \gamma }{k+1}
\\&=&\frac{\sum_{i=0}^k \int_{X_j} f(R_j^ix_j) \, d \mu_j }{k+1}\\
\\&=&\frac{\sum_{i=0}^k \int_{X_j} f(x_j) \, d \mu_j }{k+1}=\int_{X_j} f(x_j) \, d \mu_j.
\end{eqnarray*}
By a similar argument and using the fact that the cost function $c$ is invariant under each $R^k$ we obtain \[\int_{\tilde X} c(x_1,...,x_n)\, d \gamma_k=\int_{\tilde X} c(x_1,...,x_n)\, d \gamma.\]

Since the sequence $\{\gamma_k\}_{k \in \mathbb{N}}$ is tight, up to a subsequence, there exists $\bar \gamma \in \Pi (\mu_1,...,\mu_n)$ such that $\gamma_k$ converges weakly to $\bar \gamma.$
It follows from Lemma \ref{vilan} that $\int_{\tilde X} c \, d\bar \gamma \leq \liminf_{k \to \infty } c \, d \gamma_k$ from which together with the optimality of
$\gamma$ we obtain \[\int_{\tilde X} c \, d\bar \gamma=\int_{\tilde X} c \, d \gamma.\] To conclude  the proof we shall show that $R\#\bar \gamma=\bar \gamma.$
 Take a bounded continuous function $f$ on $\tilde X.$  We show that $\int_{\tilde X} f(x) \, d (R\# \bar \gamma)=\int_{\tilde X} f(x) \, d  \bar \gamma.$
Define the average of $f$ as in the Birkhoff's theorem by
\[A_{k+1}(f(x))=\frac{\sum_{i=0}^{k}f (R^i(x))}{k+1},\]
and note that
\[\int_{\tilde X}A_{k+1}(f(x)) \, d \gamma=\int_{\tilde X}f(x) \, d \gamma_k.\]
It follows from the Birkhoff's theorem that $A_{k+1}(f)$ converges $\gamma$-a.e. to an integrable function $\hat f$ and $\hat f \circ R=\hat f$ for $\gamma$-a.e. on $\tilde X.$
 Since, $f$ is bounded then so is $A_{k+1}(f)$
and therefore by the dominated convergence theorem we have $\lim_{k \to \infty}\int_{\tilde X}A_{k+1}(f(x)) \, d \gamma=\int_{\tilde X}\hat f(x) \, d \gamma.$ It implies that
\[\int_{\tilde X}\hat f(x) \, d \gamma=\lim_{k \to \infty}\int_{\tilde X}A_{k+1}(f(x)) \, d \gamma= \lim_{k \to \infty}\int_{\tilde X}f(x) \, d \gamma_k.\]
Therefore,
\begin{eqnarray*}
\int_{\tilde X} f(x) \, d (R\# \bar \gamma)=\int_{\tilde X} f(Rx) \, d \bar  \gamma
=\lim_{k \to \infty}\int_{\tilde X}f(Rx) \, d \gamma_k=\int_{\tilde X}\hat f(Rx) \, d \gamma
\end{eqnarray*}
and since $\hat f = \hat f \circ R$ $\gamma$-a.e., it follows that
\begin{eqnarray*}
\int_{\tilde X}\hat f(Rx) \, d \gamma=\int_{\tilde X}\hat f(x) \, d \gamma=\lim_{k \to \infty}\int_{\tilde X}f(x) \, d \gamma_k=\int_{\tilde X} f(x) \, d\bar \gamma.
\end{eqnarray*}
This completes the proof of the first  part. \\

Proof of part (ii): For each nonnegative integer $k$ and each $m \in\mathbb{N}$, set $U_m=(R_m,...,R_m)$ and  $U_m^k=(R_m^k,...,R_m^k).$
It follows from part (i) that there exists an optimal plan $\pi_1$ such that $U_1 \# \pi_1=\pi_1.$
 For each $k \in \mathbb{N}\cup\{0\},$ define
\[\gamma_k=\frac{\sum_{i=0}^k U_2^i\# \pi_1 }{k+1}.\]
It follows that $U_1 \# \gamma_k=\gamma_k$ for all $k \in \mathbb{N}\cup\{0\}.$ In fact, for a bounded continuous function $f$ on $\tilde X,$ we have
\begin{eqnarray*}
\int_{\tilde X} f(x) \, d (U_1 \# \gamma_k)&=&\int_{\tilde X} f(U_1 x) \, d  \gamma_k\\
&=& \frac{\sum_{i=0}^k \int_{\tilde X} f(R_1 x_1,...,R_1 x_n)\, d( U_2^i\# \pi_1) }{k+1}\\&=&
\frac{\sum_{i=0}^k \int_{\tilde X} f( R_1 R^i_2x_1,...,R_1 R^i_2x_n)\, d\pi_1 }{k+1},
\end{eqnarray*}
and since $R_2 \circ R_1=R_1 \circ R_2$, it follows that
\begin{eqnarray*}
\frac{\sum_{i=0}^k \int_{\tilde X} f( R_1 R^i_2x_1,...,R_1 R^i_2x_n)\, d\pi_1 }{k+1}&=&\frac{\sum_{i=0}^k \int_{\tilde X} f(R^i_2 R_1 x_1,...,R^i_2 R_1x_n)\, d\pi_1 }{k+1}\\
&=&\frac{\sum_{i=0}^k \int_{\tilde X} f(R^i_2  x_1,...,R^i_2 x_n)\, d\pi_1 }{k+1}=\int_{\tilde X} f(x) \, d\gamma_k.
\end{eqnarray*}
This  implies that $U_1 \# \gamma_k=\gamma_k.$  Since $\{\gamma_k\}$ is  tight, up to a subsequence,  $\gamma_k$ converges weakly to some $\pi_2$.  It then follows that $U_1\# \pi_2=\pi_2,$
 $U_2 \# \pi_2=\pi_2,$  and $\pi_2$ is a solution of $(MK).$ By repeating this argument, we obtain a sequence of  optimal transport plans $\{\pi_m\}_{m \in \mathbb{N}}$ such that
$U_k\# \pi_m=\pi_m$ for all $k \in \{1,...,m\}.$ The tightness of $\{\pi_m\}$ and Lemma \ref{vilan} ensure  the existence of  an
optimal  plan $\pi$ such that, up to a subsequence, $\pi_m$ converges weakly to $\pi.$ It also follows that $U_m \#\pi=\pi$ for every $m \in \mathbb{N}.$
This completes the proof of Theorem when there exists a finite transport plan.\\  If there is no finite transport plan then the plans $\bar \gamma$ in part (i) and $\pi$ in part (ii)
possess the desired symmetry and $\int c \, d \bar \gamma=\int c \, d \pi=\infty.$
\hfill $\square$\\

Let  $X$ be   a metric space with a metric $d$. Assume that ${\mathcal G}$ is a set of measure preserving maps on $X.$  One can equip the set $\mathcal {G}$
with the point-wise topology induced by the metric $d$. Indeed, for a sequence $\{U_k\}_{k \in \mathbb{N}} \subset {\mathcal G},$ we say that  $U_k$ converges to some $U \in \mathcal {G}$ if and only if
for every $x \in X,$
\[ \lim_{k \to \infty} d(U_kx,Ux)=0. \]
The set ${\mathcal G}$ equipped with the point-wise topology is called separable if it has a dense countable subset.\\

\textbf{Proof of Corollary \ref{mkc}.}  By assumption the cost function $c$ is invariant by a permutation $ \sigma : X^n \to X^n,$  By part (ii) of Theorem \ref{MK},
$(MK)$ has a solution $ \gamma$ such that   $ \gamma=(U,..., U)\# \gamma$ for all $U  \in \mathcal {G}.$
Define
\[\bar \gamma =\frac{\sum_{i=i}^n \sigma ^i \# \gamma}{n}.\]
It is easy to check that $\sigma \# \bar \gamma=\bar \gamma,$ and  $ \bar \gamma=(U,..., U)\#\bar \gamma$ for every  $U  \in \mathcal {G}.$\\

If $(X,d)$ is a metric space and   $\mathcal {G}$ equipped with the point-wise topology on $X$  is separable then $\mathcal {G}$ has a dense countable subset $\{U_k\}_{k \in \mathbb{N}}.$
Thus, by the latter argument,  $(MK)$ has   a solution $ \bar \gamma$ that is invariant under $\sigma,$ and  that   $ \bar \gamma=(U_k,..., U_k)\#\bar \gamma$ for
 all $ k \in \mathbb{N}.$
We shall show that $ \bar \gamma=(U,..., U)\#\bar \gamma$ for every  $U \in \mathcal {G}.$ Let $f$  be a bounded continuous function  on $X^n$, and let $U \in \mathcal {G}.$
There exists a subsequence
$\{U_{k_m}\}_{(m \in \mathbb{N})}$ approaching $U$ in the point-wise topology.
We also have
\begin{eqnarray*}
\int_{X^n} f(U_{k_m} x_1,...,U_{k_m} x_n)\, d\bar \gamma=\int_{X^n} f( x_1,..., x_n)\, d\bar \gamma.
\end{eqnarray*}
By the dominated convergence theorem we obtain
\[
\int_{X^n} f(U x_1,...,U x_n)\, d\bar \gamma=\int_{X^n} f( x_1,..., x_n)\, d\bar \gamma.
\]
and therefore $ \bar \gamma=(U,..., U)\#\bar \gamma.$ \hfill $\square$

\section{Applications}
In this section we provide two applications of our results. The first one is concerned with a volume maximizing problem.
In the second example  the cost function under the study is a repulsive function.

\subsection{Optimal transportation for the determinant.}
Given $n$ probability
measures $ \mu_1, \mu_2,..., \mu_n$ with compact support on $\R^n$,  we shall consider the following problem,
 \[\sup_{\gamma\in \Pi(\mu_1, \mu_2,..., \mu_n)} \int_{R^{n\times n}} \det(x_1,x_2,..., x_n)\, d \gamma(x_1,...,x_n). \qquad \qquad (MK_d)\]
 We focus on $(MK_d)$ in the case where the measures $(\mu_i)$ are
radially symmetric. This means that    for all $U$ in the
orthogonal group of $\R^n$,
\[U \# \mu_i=\mu_i, \qquad \forall i \in {1,...,n}.\]
Denote by $|.|$ the Euclidian norm in $\R^n$.
The problem  $(MK_d)$ is dual to :
\[\inf_{(\phi_1,...,\phi_n)\in K_d }\sum_{i=1}^n \int_{R^n} \phi_i(x_i)\, d \mu_i, \qquad (DK_d) \]
where $K_d$ is the set of n-tuples of lower-semi continuous functions $(\phi_1,...,\phi_n)$
from $\R^n$   to $\R \cup \{+ \infty\}$  such that
\[\sum_{i=1}^n  \phi_i(x_i) \geq  \det(x_1,...,x_d), \qquad \forall (x_1,...,x_n) \in \R^{n\times n}.\]
The above problem is studied in \cite{C-N}. In the radial case they  provided an explicit solution for
$(MK_d)$ from which an extremality  condition was established and  solutions of $(DK_d)$ were studied.
We shall use Theorem \ref{main} to find an explicit solution for $(DK_d)$ and derive an extremality   condition for $(MK_d).$
  Here is our result for problem $(MK_d)$ and its dual $(DK_d)$.
\begin{theorem}\label{deti} Let $\mu_1,...,\mu_n$ be radially symmetric compact supported  measures on $\R^N$ with regular  supports. If  $\mu_i$ is absolutely continuous with
respect to the  $n$-dimensional Lebesgue measure, ${\mathcal L}^n,$  and has a positive density with respect to ${\mathcal L} ^n$,
then  the following statements hold.
\begin{enumerate}
\item The dual problem has a unique solution $(\phi_1,..., \phi_n)$ such that each $\phi_i$ is radial.

\item If $\mu_1=...=\mu_n$ then we also have the following assertions:
\begin{itemize}
\item[{i.}] Set $\phi_i(x)=|x|^n/n$ for $i\in \{1,...,n\}.$ The n-tuple  $(\phi_1,...,\phi_n)$ is the unique  solution of $(DK_d).$
\item[{ii.}] Let $\bar \gamma$ be  a solution of  $(MK_d).$ For each $(x_1,..., x_n)$ in the support of $\bar \gamma,$  $(x_1,...,x_n)$ is a directed orthogonal basis in $\R^n$ with \[|x_1|=|x_2|=...=|x_n|.\]
\end{itemize}
\end{enumerate}
\end{theorem}

\textbf{Proof.} Let $B_i=\text{supp} (\mu_i).$ It follows from Proposition \ref{uniq} that $(DK_d)$ has a unique solution $(\phi_1,...,\phi_n)$
 with
 \begin{eqnarray}\phi_j(x_j)=\sup \Big \{\det(x_1, x_2,...,x_n)-\sum_{i=1, i\not=j}^n \phi_i(x_i); \, x_i \in B_i,   \, \& \,  i \not =j \Big\}.\end{eqnarray}
 The above expression shows that each $\phi_i$ is a finite convex function.
 Suppose $U$ is a  rotation matrix.  Since each $\mu_i$ is radial one has $U \# \mu_i=\mu_i.$
Note also that $\det (x_1,x_2,...,x_n)=\det(U x_1,...,U x_n)$. It follows from  Theorem \ref{main}
that $\phi_i(U x)=\phi_i(x)$ for $\mu_i$-a.e.  $x \in B_i.$ Since $\phi_i$ is continuous  and $\mu_i$ is absolutely continuous with
respect to the  $n$-dimensional Lebesgue measure, we have that  $\phi_i(U x)=\phi_i(x)$ for all   $x \in B_i.$ Thus, $\phi_i$ is invariant by 
all rotation matrices. We may now observe that,
since each  $\phi_i$ is determined by its behavior on the
half-line $\{\beta e_1; \beta \leq 0 \},$ where $e_1=(1,0,...,0) \in \R^n,$  one can write
 $\phi_i(x) = h_i(|x|),$
where the function $h_i : \R^+ \to \R$  is defined by $h_i(\beta) = \phi_i(\beta e_1).$ This competes the proof of part (1).\\

Let now  $\sigma: \R^{n\times n} \to \R^{n\times n}  $ be  a  permutation given by $\sigma(x_1, x_2,...,x_n)=(x_2,x_3,...,x_n, x_1).$
If $n$ is an odd number then $\det(x_1,...x_n)=\det \big( \sigma(x_1,...,x_n) \big).$
It now follows from part (i) of Theorem \ref{main2}, with $R$ being the identity map,   that $\phi_i=\phi_1$ for all $i \in \{1,...,n\}.$
If $n$ is an even number, let $R: \R^{n} \to \R^{n} $ be the   permutation $R(a_1, a_2,...,a_n)=(a_2,a_3,...,a_n, a_1)$ and observe that $R\# \mu_1=\mu_1 $ and 
 \[\det\big(\sigma (R x_1, R x_2,..., R x_n)\big)=\det(x_1,...,x_n).\] It now  follows from part (i) of Theorem \ref{main2} that $\phi_i(x)=\phi_1(R^{i-1}x)$ for all $x \in B.$  Thus, for both even and odd dimensions one has $\phi_i(x)=h_1(|x|)$ for $i \in \{1,...,n\}.$\\

We now show that $h_1(|x|) \geq |x|^n/n.$ For an arbitrary $x_1 \in B,$ one can find vectors $x_2,...,x_n$ so that  $(x_1,x_2,...,x_n)$ is a directed orthogonal basis of $\R^n$ and $|x_1|=...=|x_n|.$ It follows that
\[\sum_{i=1}^n h_1(|x_i|)=\sum_{i=1}^n \phi_i(x_i)\geq \det(x_1,...,x_n)=\prod_{i=1}^n|x_i|=|x_1|^n,\]
and consequently  $h_1(|x|)\geq |x|^n/n$ for all $x \in B.$
On the other hand, it follows from the Young inequality that
\[\sum_{i=1}^n \frac{|x_i|^n}{n} \geq \prod_{i=1}^n|x_i| \geq \det (x_1,...,x_n),\]
from which we obtain  $(\psi_1,...,\psi_n) \in K_d$ where $\psi_i(x)=|x|^n/n$ for every  $i \in\{1,...,n\}.$
 Since $(\phi_1,...,\phi_n)$ is a solution of $(DK_d),$ it follows that
 \[\sum_{i=1}^n \int_B\frac{|x_i|^n}{n} \,d \mu_i\geq \sum_{i=1}^n \int_B\phi_i(x_i) \,d \mu_i.\]
 The latter  inequality  together with the fact that $\phi_i(x) \geq |x|^n/n$  imply that  $\phi_i(x)=|x|^n/n$ for every  $i \in \{1,...,n\}.$ \\

 We now prove part (ii).  Since $\bar \gamma $ is a solution of $(MK_d),$ it follows from the duality relation between $(MK_d)$ and $(DK_d)$ together with part (i) of the current Theorem that
 \begin{eqnarray*}\int_{B^n}\det (x_1,...,x_n) \, d \bar \gamma&=&\sum_{i=1}^n \int_B \frac{|x_i|^n}{n} \, d \mu_i\\
&=&\sum_{i=1}^n \int_B \frac{|x_i|^n}{n} \, d \bar \gamma
  \end{eqnarray*}
  and therefore,
  \[\det(x_1,...x_n)=\sum_{i=1}^n  \frac{|x_i|^n}{n}, \qquad \bar \gamma \, \, a.e. \, \, (x_1,..., x_n) \in B^n.\]
  On the other hand by the Young inequality,
  \[\det (x_1,...x_n) \leq  \prod_{i=1}^n|x_i| \leq \sum_{i=1}^n \frac{|x_i|^n}{n}, \qquad \forall (x_1,..., x_n) \in B^n. \]
  Therefore,  on the support of $\bar \gamma $, the Young inequality becomes an equality from which part (ii) follows. \hfill $\square$

\subsection{The Coulomb cost.}
The Coulomb cost is the repulsive part of the  Hohenberg-Kohn functional and is given by \[c(x_1,...,x_n) =  \sum^n_{i\not = j} \frac{1}{|x_i-x_j|}.\]
It  represents  the Coulombic
interaction energy between the electrons. Therefore, the marginals, which represent
the single particle densities of the electrons, are all the same, embodying the indistinguishability
of the electrons.  To precisely formulate this problem, fix a probability measure $\rho$ on $\R^d.$ Let $\Pi(\rho)$ be the set of all probability measures
 on $\R^{dn}$ whose marginals are all $\rho.$ We shall consider the following problem,
 \[\inf_{\gamma\in \Pi(\rho)} \int_{R^{d\times n}} \sum^n_{i\not =j} \frac{1}{|x_i-x_j|}\, d \gamma(x_1,...,x_n). \qquad \qquad (MK_c)\]
 We focus on $(MK_c)$ in the case where the measure $\rho$ is
radially symmetric.
The problem  $(MK_c)$ is dual to :
\[\sup_{(\phi_1,...,\phi_n)\in K_c }\sum_{i=1}^n \int_{R^n} \phi_i(x_i)\, d \rho, \qquad (DK_c) \]
where $K_c$ is the set of n-tuples of lower-semi continuous functions $(\phi_1,...,\phi_n)$
from $\R^d$   to $\R \cup \{- \infty\}$  such that
\[\sum_{i=1}^n  \phi_i(x_i) \leq \sum^n_{i\not =j} \frac{1}{|x_i-x_j|}, \qquad \forall (x_1,...,x_n) \in \R^{2\times n}.\]
The optimal mass transport problem with the Coulomb cost has been recently studied  by many authors. We refer the interested reader to \cite{B-P-G, C-D-D, P2} and references therein.
Here is our result for the Coulomb cost.

\begin{theorem}\label{col} Assume that $d=2.$ Then the  following statements hold:\\
i. $(MK_c)$ has a solution $\pi \in \Pi(\rho)$ which is invariant under the rotation  group of $\R^2$, i.e., for each rotation  matrix  $U$ in $SO(2)$,
\[(U,...,U)\# \pi=\pi.\]
ii. If $\rho$ is absolutely continuous with respect to the $2$-dimensional Lebesgue measure and has a positive density function, then
 $(DK_c)$ has a solution $(\phi_1,...,\phi_n)$ where $\phi_1=...=\phi_n$ and $\phi_1$ is radially symmetric.
\end{theorem}
We refer the interested reader to \cite{G-M} where the existence of optimal transport maps is studied and a detailed proof of the above theorem is provided. Here, we shall
 sketch the proof.\\

\textbf{Proof of Theorem \ref{col}.}  Note first that for all $U , R \in  SO(2)$ one has $U \circ R=R \circ U$ as they are $2-$dimensional 
rotation matrices. Note also  that $SO(2)$ has a dense countable subset $\{R_1,R_2,...\}.$ It follows by Corollary \ref{mkc} that there exists an optimal map $\gamma$
such that $(U,...,U)\# \gamma=\gamma$ for each $U \in SO(2)$.\\

We shall now prove part (ii). It follows from Theorem \ref{main2} that $(DK_c)$ has a solution $(\phi_1,..., \phi_n)$
 such that $\phi_1=...=\phi_n$ and
\begin{equation}\label{di}
\phi_1(x_j)=\inf \big\{ c(x_1, x_2,...,x_n) - \sum_{i=1, i \not=j}^n\phi_1  (x_i);  x_i \in \R^d \, \& \,i \not=j\big \},
\end{equation}

It follows from Theorem (4) in \cite{B-P-G} that any solution of $(DK_c)$ satisfying (\ref{di}) is bounded and almost everywhere differentiable.
The same argument as in the proof of Proposition
\ref{uniq} shows that   solutions of $(DK_c)$ satisfying (\ref{di}) are unique (up to the  addition of constants summing to $0$ to each potential). It then follows from part (ii) of Theorem \ref{main2} that $\phi$ is
invariant under the group $SO(2).$ It proves that $\phi$ has to be a radial function. \hfill $\square$ \\

\end{document}